\documentclass[a4paper,12pt]{article}
\usepackage{amsmath,amsfonts,amssymb,amsthm}
\newtheorem{theo}{Th\'eor\`eme}

\newtheorem*{lem}{Lemme}
\newtheorem{conj}{Conjecture}
\begin{document}
\title{Une identit\'e en th\'eorie des partitions}
\author{Michel Lassalle\\{\small Centre National de la Recherche Scientifique}\\
{\small Ecole Polytechnique}\\
{\small 91128 Palaiseau, France}\\
{\small e-mail: lassalle @ chercheur.com}}
\date{}
\maketitle

\begin{abstract}
We prove an identity about partitions with a very elementary 
formulation. We had previously conjectured this identity, encountered in the study of 
shifted Jack polynomials. The proof given is using a trivariate generating 
function. It would be interesting to obtain a bijective proof. We 
present a conjecture generalizing this identity.
\end{abstract}

\section{Notations}

Nous revenons dans cet article sur une conjecture que nous avons 
pr\'esent\'ee  dans un pr\'ec\'edent travail (~\cite{La1}, voir aussi 
\cite{La2} ) . Il s'agit d'une identit\'e  
qui se formule de mani\`ere extr\^emement \'el\'ementaire dans 
le cadre de la th\'eorie  classique des partitions.  

Une partition $\lambda$ est une suite d\'ecroissante finie d'entiers positifs. On dit 
que le nombre
 $n$ d'entiers non nuls est la longueur de $\lambda$. On note
$\lambda  = ( {\lambda }_{1},...,{\lambda }_{n})$ 
et $n = l(\lambda)$. On dit que 
$\left|{\lambda }\right| = \sum\limits_{i = 1}^{n} {\lambda }_{i}$
est le poids de $\lambda$, et pour tout entier $i\geq1$ que 
${m}_{i} (\lambda)  = card \{j: {\lambda }_{j}  = i\}$
est la multiplicit\'e de $i$ dans $\lambda$. On identifie $\lambda$ \`a son diagramme 
de Ferrers 
$\{ (i,j) : 1 \le i \ \le l(\lambda), 1 \ \le j \ \le {\lambda }_{i} \}$. 
On pose 
\[{z}_{\lambda }  = \prod\limits_{i \ge  1}^{} {i}^{{m}_{i}(\lambda)} {m}_{i}(\lambda) 
!  .\]
	 
Nous avons introduit dans~\cite{La1} la g\'en\'eralisation suivante du coefficient 
binomial classique. Soient $\lambda$ une partition et $r$ un entier $\geq 1$. On note 
$\genfrac{\langle}{\rangle}{0pt}{}{\lambda}{r}$ 
le nombre de fa\c cons dont on peut choisir r points 
dans le diagramme de $\lambda$ de telle sorte que 
\textit{au moins un point soit choisi sur chaque ligne de $\lambda$}.
	
Les coefficients binomiaux g\'en\'eralis\'es  
$\genfrac{\langle}{\rangle}{0pt}{}{\lambda}{r}$ poss\`edent la fonction 
g\'en\'eratrice suivante
\begin{equation*}
\begin{split}
\sum_{r \ge  1} \genfrac{\langle}{\rangle}{0pt}{}{\lambda}{r}\,{q}^{r} &= 
\prod_{i = 1}^{l(\lambda)} \left({{(1 + q)}^{{\lambda}_{i}} - 1}\right)\\
 &= \prod_{i \ge  1}^{} {\left({{(1 + q)}^{i} - 1}\right)}^{{{m}_{i}(\lambda)}}.
\end{split}
\end{equation*}

Soient $X$ une ind\'etermin\'ee et $n$ un entier $\geq1$. On note d\'esormais 
\[{(X)}_{n }  =  X (X +1) ... (X + n-1)\] 
\[{[X]}_{n }  =  X (X - 1) ... (X -n +1)  .\]	
les factorielles ``ascendante'' et ``descendante'' classiques. On pose 
\[\binom{X}{n}  =  {\frac{{[X]}_{n}}{n!}}  .\]
Les nombres de Stirling de premi\`ere esp\`ece $s(n,k)$ sont d\'efinis 
par les fonctions g\'en\'eratrices
\[{[X]}_{n} =\sum_{k \ge  1} s(n,k){X}^{k}\]
\[{(X)}_{n} =\sum_{k \ge  1} \mid s(n,k) \mid {X}^{k}.\]

\section{Notre r\'esultat}

Il s'agit d'une g\'en\'eralisation de la propri\'et\'e classique suivante, qui est par 
exemple d\'emontr\'ee au Chapitre 1, Section 2, Exemple 1 du livre de 
Macdonald~\cite{Ma}.
Soit $X$ une ind\'etermin\'ee. Pour tout entier  $n\geq1$ on a 
\[\sum_{\left|{\mu }\right| = n} {(-1)}^{n-l(\mu)} {\frac{{X}^{l(\mu)}}{{z}_{\mu }}} = \binom{X}{n}\]  
\[\sum_{\left|{\mu }\right| = n} {\frac{{X}^{l(\mu)}}{{z}_{\mu }}}  = 
\binom{X + n - 1}{n}.\]	
Ces deux relations sont \'equivalentes en changeant $X$ en $-X$. 
	
Dans~\cite{La1} nous avons formul\'e la conjecture 
suivante: pour tous entiers $n,r,s\geq1$ on a
\begin{multline*}
\sum_{\left|{\mu }\right| = n} {(-1)}^{r - l(\mu )} 
{\frac{\displaystyle\genfrac{\langle}{\rangle}{0pt}{}{\mu}{r}}
{{z}_{\mu }}} {X}^{l(\mu ) - 1} \left({ \sum_{i = 1}^{l(\mu )} 
{({\mu }_{i})}_{s} }\right)  \\
= (s - 1) ! \binom{n + s - 1}{n - r} \sum_{i = 1}^{\min (r,s)}  
\binom{X - s}{r - i} \binom{s}{i} .
\end{multline*}
Comme on a
$\genfrac{\langle}{\rangle}{0pt}{}{\mu}{\left|{\mu }\right|}  = 1$,
on voit facilement que la propri\'et\'e classique pr\'ec\'edente 
correspond au cas $r=n$, $s=1$.

Compte-tenu de la formule classique de Chu-Vandermonde
\begin{equation}
\binom{X+Y}{n}= \sum_{i = 0}^{n}\binom{X}{n-i}\binom{Y}{i}
\end{equation}
notre conjecture peut s'exprimer sous la forme plus naturelle suivante.
\begin{theo} 
Soit $X$ une ind\'etermin\'ee. Pour tous entiers $n,r,s\geq1$ on a
\begin{multline*}
\sum_{\left|{\mu }\right| = n} {(-1)}^{r - l(\mu )} 
{\frac{\displaystyle\genfrac{\langle}{\rangle}{0pt}{}{\mu}{r}}{{z}_{\mu }}} 
{X}^{l(\mu ) - 1} \left({ \sum_{i = 1}^{l(\mu )} {({\mu }_{i})}_{s} }\right) \\
= (s - 1) ! \binom{n + s - 1}{n - r} \left[{\binom{X}{r} - \binom{X - s}{r}}\right] .
\end{multline*}	
Ou de mani\`ere \'equivalente, en  changeant $X$  en $-X$,
\begin{multline*}
\sum_{\left|{\mu }\right| = n} 
{\frac{\displaystyle\genfrac{\langle}{\rangle}{0pt}{}{\mu}{r}}{{z}_{\mu }}} 
{X}^{l(\mu ) - 1} \left({ \sum_{i = 1}^{l(\mu )} {({\mu }_{i})}_{s} }\right) \\
= (s - 1) ! \binom{n + s - 1}{n - r} 
\left[{\binom{X + r + s - 1}{r} - \binom{X + r - 1}{r}}\right] .
\end{multline*}
\end{theo}

Comme on a
\[\genfrac{\langle}{\rangle}{0pt}{}{\mu}{r}  = 0   
\qquad\textrm{si } r > \left|{\mu }\right| ,\]
l'identit\'e est triviale pour $r>n$. Comme on a
\[\genfrac{\langle}{\rangle}{0pt}{}{\mu}{r}  = 0  \qquad\textrm{si } r  < l(\mu ) ,\]
la sommation au membre de gauche est limit\'ee aux partitions 
$\mu$ telles que $l(\mu)\leq r$.
Chacun des membres de l'identit\'e est ainsi un polyn\^ome en $X$ de degr\'e  $r-1$.
	
Le Th\'eor\`eme 1 se d\'ecompose donc en $r$ identit\'es obtenues en 
identifiant les coefficients de ${X}^{p} ( 0  \le  p  \le  r - 1 )$ dans chaque 
membre. Le coefficient de ${X}^{p-1}$ est obtenu par sommation sur les partitions de 
longueur $p$. De mani\`ere \'equivalente le Th\'eor\`eme 1 se formule donc comme suit.
\begin{theo}
Pour tous entiers $n,r,s\geq1$ et pour tout entier $1 \le  p  \le  r $ on a
\begin{multline*}
r !  \sum_{\begin{subarray}{1}\left|{\mu }\right| = n\\l(\mu ) = p\end{subarray}} 
{\frac{\displaystyle\genfrac{\langle}{\rangle}{0pt}{}{\mu}{r}}{{z}_{\mu }}} 
{ \left( \sum_{i = 1}^{l(\mu )} \frac{{({\mu }_{i})}_{s}}{s!}\right)} \\ 
=  \binom{n + s - 1}{n - r} 
\left(\sum_{j=p}^{r} \binom{j}{p-1} \mid s(r,j) \mid {s}^{j-p}\right).
\end{multline*}
\end{theo}

Dans le cas particulier $s=1$, le Th\'eor\`eme 1 s'\'ecrit facilement (voir \cite{La1}, p. 462)
\begin{equation}
\sum_{\left|{\mu }\right| = n} {(-1)}^{r-l(\mu)} 
\frac{\displaystyle{\genfrac{\langle}{\rangle}{0pt}{}{\mu}{r}}}{z_{\mu }} X^{l(\mu)}
= \binom{n-1}{r-1} \binom{X}{r}.
\end{equation}
Le Th\'eor\`eme 2 devient alors
\begin{equation}
r !  \sum_{\begin{subarray}{1}\left|{\mu }\right| = n\\l(\mu ) = p\end{subarray}} 
{\frac{\displaystyle\genfrac{\langle}{\rangle}{0pt}{}{\mu}{r}}{{z}_{\mu }}} 
=  \binom{n - 1}{r - 1} \mid s(r,p) \mid .
\end{equation}

Rodica Simion a obtenu une preuve purement combinatoire de cette 
identit\'e  \cite{Ro}. 
Sa d\'emonstration est donn\'ee dans \cite{La1}, p. 461. Il s'agit d'une preuve 
bijective ``\`a la Sch\"utzenberger'': on montre que chaque membre de (3) 
d\'enombre un certain ensemble de 
permutations de deux fa\c  cons diff\'erentes. Il serait int\'eressant de 
disposer d'une preuve de ce type dans le cas g\'en\'eral $s\neq 1$. 

Cet article n'aurait pas pu \^etre \'ecrit sans les remarques 
et les conseils d'Alain Lascoux, qui nous a montr\'e comment traduire 
notre conjecture en termes de fonctions sym\'etriques, ce qui est mis en 
oeuvre dans les sections 4 et 5. Nous pouvons ainsi donner deux preuves 
 diff\'erentes du Th\'eor\`eme 2 aux sections 6 et 7.
  
Alain Lascoux nous a \'egalement signal\'e le lien \'etroit entre notre identit\'e
et la formule de Cauchy, ainsi que sa formulation en termes de $\lambda$ -anneaux, 
 qui est esquiss\'ee \`a la section 8.
 
Enfin la Section 9 expose bri\`evement les motivations nous ayant 
 conduit au Th\'eor\`eme 1. Il s'agit de l'\'etude des polyn\^omes ``sym\'etriques 
 d\'ecal\'es''. Nous formulons en particulier une nouvelle conjecture qui g\'en\'eralise 
 l'identit\'e du Th\'eor\`eme 1. 
 
\section{Deux cas particuliers}

Il est tr\`es facile d'expliciter le Th\'eor\`eme 2   pour $p=r$, et $p=r-1$. 
Pour $p=r$ le Th\'eor\`eme 2 s'\'ecrit
\[(r-1)!\sum_{\begin{subarray}{1}\left|{\mu }\right| = n\\l(\mu ) = r\end{subarray}} 
{\frac{\displaystyle\genfrac{\langle}{\rangle}{0pt}{}{\mu}{r}}{{z}_{\mu }}} 
{ \left({ \sum_{i = 1}^{l(\mu )} ({\mu }_{i})_{s} }\right)}
= s! \binom{n + s - 1}{n - r}.\] 
Comme  on a 
\[\genfrac{\langle}{\rangle}{0pt}{}{\mu}{l(\mu )}  = \prod_{i = 1}^{l(\mu )} 
{\mu }_{i} = \prod_{i \ge 1} {i}^{{m}_{i}(\mu )} ,\]	
on obtient le
\begin{theo} 
Pour tous entiers $n,r,s\geq1$  on a
\[(r - 1) !  \sum_{\begin{subarray}{1}\left|{\mu }\right| = n\\
l(\mu ) = r\end{subarray}} {\frac{\displaystyle\sum\limits_{i \ge  1}^{} {{m}_{i}(\mu ) (i)}_{s} }{\displaystyle\prod\limits_{i \ge  1}^{} {m}_{i}(\mu ) !}} 
=  s ! \binom{n + s - 1}{n - r} .\]
\end{theo}

Pour $p=r-1$ le Th\'eor\`eme 2 s'\'ecrit 
\[2(r-2)!\sum_{\begin{subarray}{1}\left|{\mu }\right| = n\\l(\mu ) = r-1\end{subarray}} 
{\frac{\displaystyle\genfrac{\langle}{\rangle}{0pt}{}{\mu}{r}}{{z}_{\mu }}} 
{ \left({ \sum_{i = 1}^{l(\mu )} ({\mu }_{i})_{s} }\right)}
= s!\binom{n + s - 1}{n - r} (s+r-1)\] 
car on a $s(r,r-1)=-\binom{r}{2}$. On v\'erifie facilement la relation 
\[\genfrac{\langle}{\rangle}{0pt}{}{\mu}{l(\mu ) + 1} =  
{\frac{1}{2}} \left({\left|{\mu }\right| - l(\mu )}\right)
 \prod_{i = 1}^{l(\mu )} {\mu }_{i} .\]	
L'identit\'e pr\'ec\'edente s'\'ecrit donc
\[(r - 2) !  \sum_{\begin{subarray}{1}\left|{\mu }\right| = n\\
l(\mu ) = r-1\end{subarray}} 
{\frac{\displaystyle\sum\limits_{i \ge  1}^{} {{m}_{i}(\mu ) (i)}_{s} }
{\displaystyle\prod\limits_{i \ge  1}^{} {m}_{i}(\mu ) !}} 
=  s ! \binom{n + s - 1}{n - r}\frac{s+r-1}{n-r+1} \]
c'est-\`a-dire la relation du Th\'eor\`eme 3, \'ecrite en rempla\c cant $r$ par $r-1$.

Enfin comme on a
\[\genfrac{\langle}{\rangle}{0pt}{}{\mu}{\left|{\mu }\right|}  = 1 \]
le Th\'eor\`eme 1 prend la forme suivante pour $r=n$.

\begin{theo} 
Soit $X$ une ind\'etermin\'ee. Pour tous entiers $n,s\geq1$ on a
\[\sum_{\left|{\mu }\right| = n} {(-1)}^{n - l(\mu )} 
{\frac{{X}^{l(\mu ) - 1}}{{z}_{\mu }}}  \left({ \sum_{i = 1}^{l(\mu )} 
{({\mu }_{i})}_{s} }\right)  =  (s - 1) !  
\left[{\binom{X}{n}  -  \binom{X - s}{n}}\right]  .\]
Ou de mani\`ere \'equivalente	
\[\sum_{\left|{\mu }\right| = n}  {\frac{{X}^{l(\mu ) - 1}}{{z}_{\mu }}}  
\left({ \sum_{i = 1}^{l(\mu )} 
{({\mu }_{i})}_{s} }\right)  =  (s - 1) !  
\left[{\binom{X + n + s - 1}{n} - \binom{X + n - 1}{n}}\right] .\]
\end{theo}

Pour $s=1$ on retrouve la propri\'et\'e classique \'enonc\'ee au 
commencement de la section 2.

\section{Le cas particulier $r=n$}

On consid\`ere les fonctions sym\'etriques compl\`etes ${h}_{i}(x)$ d'un 
ensemble de variables $x$ (~\cite{Ma}, Chapitre 1, Section 2). Rappelons que 
${h}_{i}(x)$ est le coefficient de $t^{i}$ dans le d\'eveloppement en s\'erie 
enti\`ere de $\prod_{j} {(1 - {x}_{j}t)}^{- 1}$. 
Pour deux ensembles de variables $z,u$ on a donc imm\'ediatement
\begin{equation}
{h}_{n}(z,u) =\sum_{i=0}^{n} {h}_{n-i}(z) {h}_{i}(u). 
\end{equation}
On note $(1^{k})$ le point 
${x}_{1} = \cdots = {x}_{k} = 1 , {x}_{k + 1} = {x }_{k + 2}= \cdots = 0$.
Alors on a (~\cite{Ma}, Exercice 1.1.2)
\begin {equation}
{h}_{i}(1^{k}) =\binom{i+k-1}{i}.
\end{equation}
On note ${h}_{i}(1^{X})$ le polyn\^ome obtenu par la continuation analytique du second 
membre. C'est le coefficient de $t^{i}$ dans le d\'eveloppement en s\'erie enti\`ere de
${(1 - t)}^{- X}$.

A titre d'explicitation de la m\'ethode suivie, nous 
d\'emontrons maintenant le Th\'eor\`eme 4. La seconde identit\'e peut s'\'ecrire
\[\sum_{\left|{\mu }\right| = n} {\frac{{X}^{l(\mu ) - 1}}{{z}_{\mu }}} 
 \left( \sum_{i = 1}^{l(\mu )} h_{{\mu }_{i}-1}(1^{s+1}) \right) 
= \frac{1}{s} \Big(h_{n}(1^{X+s})-h_{n}(1^{X})\Big) .\]
En appliquant (4) et la relation \'el\'ementaire
\[\frac{{h}_{i}(1^{k})}{k} =\frac{{h}_{i-1}(1^{k+1})}{i}\]
le membre de droite devient
\begin{equation*}
\frac{1}{s}  \Big(h_{n}(1^{X+s})-h_{n}(1^{X})\Big) =
\sum_{i = 1}^{n} {h}_{n - i}({1}^{X}){\frac{{h}_{i}({1}^{s})}{s}}=
\sum_{i = 1}^{n} {h}_{n - i}({1}^{X}){\frac{{h}_{i-1}({1}^{s+1})}{i}}.
\end{equation*}
Le Th\'eor\`eme 4 se formule donc comme suit
\[\sum_{\left|{\mu }\right| = n} {\frac{{X}^{l(\mu ) - 1}}{{z}_{\mu }}} 
 \left({ \sum_{i \ge 1} {m }_{i}(\mu) h_{i-1}(1^{s+1})}\right) 
= \sum_{i = 1}^{n} {h}_{n - i}({1}^{X}){\frac{{h}_{i-1}(1^{s+1})}{i}}.\]

Cette relation est le cas 
particulier pour $z=1^{s+1}$ de l'identit\'e \textit{plus g\'en\'erale} 
suivante, \'ecrite pour un ensemble \textit{quelconque} de variables $z$,
\[\sum_{\left|{\mu }\right| = n} {\frac{{X}^{l(\mu ) - 1}}{{z}_{\mu }}} 
 \left({ \sum_{i \ge 1} {m }_{i}(\mu) h_{i-1}(z)}\right) 
= \sum_{i = 1}^{n} {h}_{n - i}({1}^{X}){\frac{{h}_{i-1}(z)}{i}}.\]

Maintenant, et c'est le point \textit{essentiel}, on peut observer que cette nouvelle 
identit\'e est \textit{lin\'eaire} en les fonctions ${h}_{i}(z)$. 
Il suffit donc de la d\'emontrer lorsque $z$ est r\'eduit \`a \textit{une seule 
variable}. On a alors ${h}_{i}(z) = {z}^{i}$.  

L'identit\'e pr\'ec\'edente est ainsi \textit{\'equivalente} \`a celle dont nous donnons 
maintenant la d\'emonstration.

\begin{theo} 
Soient $X$ et $z$ deux ind\'etermin\'ees. Pour tout entier $n \ge 1$ on a
\[\sum_{\left|{\mu }\right| = n} {\frac{{X}^{l(\mu ) - 1}}{{z}_{\mu }}} 
 \left( \sum_{i \ge  1} {m}_{i}(\mu ) {z}^{i-1} \right) =
\sum_{i = 1}^{n} {h}_{n - i}({1}^{X}) \frac{{z}^{i-1}}{i}.\]	
\end{theo}

\begin{proof}[Preuve]
On forme la fonction g\'en\'eratrice
\[\sum_{\mu } {X}^{l(\mu )}{y}^{\left|{\mu }\right|}\            
{\frac{\displaystyle\sum_{i \ge  1} {{m}_{i}(\mu ) {z}^{i-1}} }{\displaystyle\prod_{i \ge  1} {i}^{{m}_{i}(\mu )}{m}_{i}(\mu ) !}}.\]
On peut l'\'ecrire comme la d\'eriv\'ee, prise au point $u=1$ de
\[\sum_{m_i} {{\frac{1}{\prod_{i} {i}^{{m}_{i}}{m}_{i} !}}
\,{X}^{\sum_{i} 
{m}_{i}}\,{y}^{\sum_{i} i{m}_{i}}u^{\sum_{i} {{m}_{i} {z}^{i-1}}}}.\]
La fonction g\'en\'eratrice est donc \'egale \`a 
\begin{equation*}
\begin{split}
{d\over du}  \sum_{m_i} \prod_{i}  \frac {1}{m_i !}{\bigg( X\,\frac{y^i}{i}\,u^{z^{i-1}}\bigg)}^{m_i} \ \Big|_{u=1}
&= {d\over du} \prod_{i} exp( X\, \frac{y^i}{i}\, u^{z^{i-1}})\Big|_{u=1}\\
&= \left( X\sum_{i\ge 1} z^{i-1}  \frac{y^i}{i} \right) exp\big(X\,\sum_{i\ge 1} 
\frac{y^i}{i}\big)\\
&= \left( X\sum_{i\ge 1} z^{i-1}  \frac{y^i}{i} \right) \frac{1}{{(1-y)}^{X}}.
\end{split}
\end{equation*} 
On voit ainsi que le terme $y^{n}$ apparait dans chacune
 des contributions suivantes pour tout $i\ge 1$,
\[\frac{X}{z}\frac{{(zy)}^{i}}{i}\quad\textrm{x coefficient de}\quad{y^{n-i}} 
\quad\textrm{dans}\quad{(1-y)^{-X}}.\]
Par la d\'efinition m\^eme des fonctions $h_i$, le coefficient de $y^{n}$ est donc
\[X\sum_{i\ge 1} \frac{{z}^{i-1}}{i}{h}_{n-i}(1^{X}).\]
\end{proof}

\section{La m\'ethode g\'en\'erale}

Nous revenons maintenant \`a la formulation g\'en\'erale du Th\'eor\`eme 1, 
avec $r$ arbitraire. Nous aurons besoin du r\'esultat auxiliaire 
suivant.
\begin{lem}
Pour tous entiers $a,b,c,d$ on a
\[\binom{a+b-1}{d-c}\binom{b+c-1}{c-1}=\sum_{i=c}^{d}
\binom{i-1}{c-1}\binom{a-i-1}{a-d-1}
\binom{b+i-1}{i-1}.\]
\end{lem}
\begin{proof}[Preuve]
L'expression se r\'eduit imm\'ediatement \`a
\[\binom{a+b-1}{d-c}=\sum_{i=c}^{d}\binom{a-i-1}{a-d-1}
\binom{b+i-1}{b+c-1}.\]
Ce qui peut s'\'ecrire en posant $k=d-c$,
\[\binom{a+b-1}{k}=\sum_{i=0}^{k}\binom{a-c-i-1}{k-i}\binom{b+c+i-1}{i}.\]
On termine en appliquant (1).
\end{proof}

Compte-tenu de (1) on peut \'ecrire
\[\binom{X-s}{r}-\binom{X}{r}=\sum_{j=1}^{r}\binom{X}{r-j}\binom{-s}{j}.\]
Mais on a
\[\frac{1}{s}\binom{-s}{j}=\frac{{(-1)}^{j}}{s}\binom{s+j-1}{j}=\frac{{(-1)}^{j}}{j}
\binom{s+j-1}{j-1}.\]
Compte-tenu de (5) l'identit\'e du Th\'eor\`eme 1 peut donc se formuler comme suit  
\begin{multline*}
\sum_{\left|{\mu }\right| = n} {(-1)}^{r-l(\mu )} {\frac{\displaystyle\genfrac{\langle}{\rangle}{0pt}{}{\mu}{r}}{{z}_{\mu }}} 
{X}^{l(\mu ) - 1} \left({ \sum_{i=1}^{l(\mu)}  h_{{{\mu}_{i}}-1}(1^{s+1})}\right)\\
= \sum_{j = 1}^{r} \frac{(-1)^{j-1}}{j}\binom{X}{r-j}\,\binom{n+s-1}{n-r}
\binom{s+j-1}{j-1}.
\end{multline*}

Maintenant, en appliquant le Lemme pour $n,s,j,n-r+j$ , on a
\[\binom{n+s-1}{n-r}\binom{s+j-1}{j-1}=\sum_{i= j}^{n-r+j}\binom{i-1}{j-1}
\binom{n-i-1}{r-j-1}\binom{s+i-1}{i-1}.\]
Le Th\'eor\`eme 1 peut donc s'\'ecrire
\begin{multline*}
\sum_{\left|{\mu }\right| = n} {(-1)}^{r-l(\mu )}
{\frac{\displaystyle\genfrac{\langle}{\rangle}{0pt}{}{\mu}{r}}{{z}_{\mu }}} 
{X}^{l(\mu ) - 1} \left({ \sum_{i\ge 1} m_{i}(\mu) h_{i-1}(1^{s+1})}\right)\\
 =\sum_{j = 1}^{r}\sum_{i= j}^{n-r+j} (-1)^{j-1}\binom{X}{r-j}\,
{\frac{1}{i}}\binom{i}{j}\binom{n-i-1}{r-j-1}\,{h}_{i-1}(1^{s+1}).
\end{multline*}

Comme pr\'ec\'edemment cette identit\'e est le cas particulier pour $z=1^{s+1}$ de 
l'identit\'e plus g\'en\'erale suivante, \'ecrite pour un ensemble 
\textit{quelconque} de variables $z$.

\begin{theo} 
Soient $X$ une ind\'etermin\'ee et $z$ un ensemble quelconque de variables. 
Pour tous entiers $n,r\geq1$ on a
\begin{multline*}
\sum_{\left|{\mu }\right| = n} {(-1)}^{r-l(\mu )}
 {\frac{\displaystyle\genfrac{\langle}{\rangle}{0pt}{}{\mu}{r}}{{z}_{\mu }}} 
{X}^{l(\mu ) - 1} \left({ \sum_{i\ge 1} m_{i}(\mu) h_{i-1}(z)}\right)\\
 =\sum_{j = 1}^{r}\sum_{i= j}^{n-r+j} (-1)^{j-1}\binom{X}{r-j}\,
{\frac{1}{i}}\binom{i}{j}\binom{n-i-1}{r-j-1}\,{h}_{i-1}(z).
\end{multline*}
\end{theo}

Comme pr\'ec\'edemment cette nouvelle identit\'e est lin\'eaire en les fonctions 
${h}_{i}(z)$. Il suffit donc de la d\'emontrer lorsque $z$ est r\'eduit 
\`a \textit{une seule variable}. On a alors ${h}_{i}(z) = {z}^{i}$.
 
Le Th\'eor\`eme 6 est ainsi \textit{\'equivalent} \`a celui que nous formulons 
maintenant.

\begin{theo} 
Soit $X$ et $z$ deux ind\'etermin\'ees. Pour tous entiers $n,r\geq1$ on a
\begin{multline*}
\sum_{\left|{\mu }\right| = n} {(-1)}^{r-l(\mu )}
 {\frac{\displaystyle\genfrac{\langle}{\rangle}{0pt}{}{\mu}{r}}{{z}_{\mu }}} 
{X}^{l(\mu ) - 1} \left({ \sum_{i\ge 1} m_{i}(\mu) z^{i-1}}\right)\\
 =\sum_{j = 1}^{r}\sum_{i= j}^{n-r+j} (-1)^{j-1}\binom{X}{r-j}\,
{\frac{1}{i}}\binom{i}{j}\binom{n-i-1}{r-j-1}\,{z}^{i-1}.
\end{multline*}
\end{theo}  

On peut formuler ce r\'esultat de mani\`ere \'equivalente en identifiant dans chaque membre 
les coefficients de $X^{p-1}$, avec $1 \le  p  \le  r $.
\begin{theo} 
Pour tous entiers $n,r\geq1$, pour tout entier $1 \le  p  \le  r $ et pour toute
ind\'etermin\'ee  $z$ on a
\begin{multline*}
\sum_{\begin{subarray}{1}\left|{\mu }\right| = n\\l(\mu ) = p\end{subarray}} 
 {\frac{\displaystyle\genfrac{\langle}{\rangle}{0pt}{}{\mu}{r}}{{z}_{\mu }}} 
 \left({ \sum_{i\ge 1} m_{i}(\mu) z^{i-1}}\right)\\
 =\sum_{j = 1}^{r}\sum_{i= j}^{n-r+j} 
{\frac{\mid s(r-j,{p-1}) \mid}{i\,(r-j)!}}\binom{i}{j}\binom{n-i-1}{r-j-1}\,{z}^{i-1}.
\end{multline*}
\end{theo}

Cet \'enonc\'e est ainsi une g\'en\'eralisation du Th\'eor\`eme 2, et c'est 
lui que nous allons maintenant d\'emontrer. 

\section{Premi\`ere d\'emonstration}

Nous aurons besoin du r\'esultat auxiliaire suivant, d\'ej\`a utilis\'e 
par Di Bucchianico et Loeb~\cite{Buc}.

\begin{theo}
Soient $x,y,q$ trois ind\'etermin\'ees. On a
\[{\left(\frac{1-y}{1-y(1+q)}\right)}^{x}=\sum_{i,j,k\ge 0}{\binom{i-1}{j-1}
\frac{\mid s(j,k) \mid}{j!}x^{k}y^{i}q^{j}}\] 
\end{theo}
\begin{proof}[Preuve]
Nous appliquons la formule du bin\^ome
\[\sum_{i\ge j}{\binom{i-1}{j-1}y^{i-j}}=\frac{1}{(1-y)^{j}}.\]
Le membre de droite s\'ecrit donc comme suit
\begin{equation*}
\begin{split}
\sum_{i,j}{\binom{i-1}{j-1}\binom{-x}{j}\,y^{i}{(-q)}^{j}}&=
\sum_{j}{\binom{-x}{j}\,{\left(\frac{qy}{y-1}\right)}^j}
\end{split}
\end{equation*}
L'\'enonc\'e r\'esulte alors de la formule du bin\^ome
\[{(1+t)}^{-x}=\sum_{j\ge 0}{\binom{-x}{j}t^j}.\]
\end{proof}

\begin{proof}[D\'emonstration du Th\'eor\`eme 8]
Nous formons la fonction g\'en\'eratrice
\[\sum_{\mu,r} {q}^{r}{x}^{l(\mu )}{y}^{\left|{\mu }\right|}            
\ {\frac{\displaystyle\genfrac{\langle}{\rangle}{0pt}{}{\mu}{r}}
{\displaystyle\prod_{i \ge  1} {i}^{{m}_{i}(\mu )}{m}_{i}(\mu ) !}}\ \sum_{i \ge  1}
 {{m}_{i}(\mu ) {z}^{i-1}}.\]
On peut l'\'ecrire comme la d\'eriv\'ee, prise au point $u=1$ de
\begin{multline*}
\sum_{m_i} u^{\sum_{i} {{m}_{i} {z}^{i-1}}}{x}^{\sum_{i} {m}_{i}}\,{y}^{\sum_{i} 
i{m}_{i}}
\prod_{i} \frac{\displaystyle {((1 + q)}^{i} - 1)^{m_i}}
{\displaystyle {i}^{{m}_{i}}{m}_{i} !}=\\
\sum_{m_i} \prod_{i}  \frac {1}{m_i !} {\bigg( x\,y^i\,u^{z^{i-1}}\,
\frac{{(1 + q)}^{i} - 1}{i}\bigg)}^{m_i}.
\end{multline*}
La fonction g\'en\'eratrice est donc \'egale \`a
\begin{multline*}
{d\over du} \prod_{i} exp\Big( x\, y^i\, \frac{{(1 + q)}^{i} - 1}{i} \, u^{z^{i-1}} \Big) \Big|_{u=1}
= \\
x\,\left( \sum_{i\ge 1} z^{i-1} y^i\,\frac{{(1 + q)}^{i} - 1}{i} \right) 
exp\left(x\,\sum_{i\ge 1}
 y^i\,\frac{{(1 + q)}^{i} - 1}{i}\right).
\end{multline*}
Mais le dernier terme peut s'\'ecrire
\begin{equation*}
\begin{split}
exp\left(x\,\sum_{i\ge 1} y^i\,\frac{{(1 + q)}^{i} - 1}{i}\right)&=
exp\left(x\,\sum_{i\ge 1} \big(\frac{{(y(1 + q))}^{i}}{i}-\frac{{y}^{i}}{i}\big)\right)\\
&=exp(x[-log(1-y(1+q)+log(1-y)])\\
&={\left(\frac{1-y}{1-y(1+q)}\right)}^{x}.
\end{split}
\end{equation*}
La fonction g\'en\'eratrice s'\'ecrit donc
\[\frac{x}{z} \left( \sum_{i\ge 1} \frac{{(yz)}^{i}}{i} 
\,\left( \sum_{j=1}^{i}\binom{i}{j}{q}^j \right) \right)
{\left(\frac{1-y}{1-y(1+q)}\right)}^{x}\]
Le Th\'eor\`eme 9 permet d'expliciter le terme en $q^rx^py^n$ de 
cette fonction. On voit que le terme $q^rx^py^n$ apparait dans chacune
des contributions suivantes pour tout $i\ge 1$ et pour tout $1\le j \le i$,
\[\frac{x}{z}\frac{{(yz)}^{i}}{i}\quad \binom{i}{j}\,{q}^j\quad
{\binom{n-i-1}{r-j-1}\frac{\mid s(r-j,p-1) \mid}{(r-j)!}{x}^{p-1}y^{n-i}{q}^{r-j}}.\]
Le coefficient de $q^rx^py^n$ est donc
\[\sum_{i=1}^{n-1} \sum_{j=1}^{i}  \binom{i}{j}\,
\frac{\mid s(r-j,p-1) \mid}{i\,(r-j)!}\,\binom{n-i-1}{r-j-1}\,z^{i-1}.\]
Ce qui ach\`eve la preuve du Th\'eor\`eme 8.
\end{proof}

\section{Seconde d\'emonstration}

Cette seconde d\'emonstration est due \`a Theresia Eisenk\"olbl 
\cite{Th2}, qui a trouv\'e une preuve tr\`es rapide du Th\'eor\`eme 8.
Fixons $i \ge 1$ et identifions le coefficient 
de ${z}^{i-1}$ dans chaque membre. Nous obtenons
\[\sum_{\begin{subarray}{1}\left|{\mu }\right| = n\\l(\mu ) = p\end{subarray}} 
 {\frac{\displaystyle\genfrac{\langle}{\rangle}{0pt}{}{\mu}{r}}{{z}_{\mu }}} 
  m_{i}(\mu) =\sum_{j = 1}^{i} 
{\frac{\mid s(r-j,{p-1}) \mid}{i\,(r-j)!}}\binom{i}{j}\binom{n-i-1}{r-j-1}.\]

La sommation au membre de gauche est r\'eduite aux partitions $\mu$ telles que $m_{i}(\mu) 
\neq 0$, c'est-\`a-dire de la forme $\mu = \lambda \cup \{i\}$. On a alors $m_{i}(\mu) 
=m_{i}(\lambda)+1$, $l(\mu )=l(\lambda )+1$, $z_{\mu}=i m_{i}(\mu) z_{\lambda}$ 
et
\[\genfrac{\langle}{\rangle}{0pt}{}{\mu}{r}= \sum_{j = 1}^{i} 
\binom{i}{j} \genfrac{\langle}{\rangle}{0pt}{}{\lambda}{r-j}.\] 
La relation du Th\'eor\`eme 8 s'\'ecrit donc
\[\sum_{\begin{subarray}{1}\left|{\lambda }\right| = n-i\\l(\lambda ) = 
p-1\end{subarray}} \sum_{j = 1}^{i} \binom{i}{j}
 {\frac{\displaystyle\genfrac{\langle}{\rangle}{0pt}{}{\lambda}{r-j}}{{z}_{\lambda }}} 
 =\sum_{j = 1}^{i} 
{\frac{\mid s(r-j,{p-1}) \mid}{(r-j)!}}\binom{i}{j}\binom{n-i-1}{r-j-1}.\]
Mais c'est exactement la relation (3) 
\[(r-j) !  \sum_{\begin{subarray}{1}\left|{\lambda }\right| = n-i\\l(\mu ) = 
p-1\end{subarray}} 
{\frac{\displaystyle\genfrac{\langle}{\rangle}{0pt}{}{\lambda}{r-j}}{{z}_{\lambda }}} 
=  \binom{n -i- 1}{r -j- 1} \mid s(r-j,p-1) \mid \]
multipli\'ee par $\binom{i}{j}$ et somm\'ee sur $j \ge 1$.\qed

Theresia Eisenk\"olbl a donn\'e dans \cite{Th1} une d\'emonstration diff\'erente 
du Th\'eor\`eme 1.

\section{Fonctions sym\'etriques et $\lambda$-anneaux}

Dans cette section nous allons consid\'erer les fonctions sym\'etriques comme 
\'etant des op\'erateurs sur l'anneau des polyn\^omes d'un nombre quelconque 
d'ind\'etermin\'ees 
$(v_{1},v_{2},\ldots,v_{n},y,t,\ldots)$. Voir par exemple \cite{Ma}, Remarque
1.2.15, et \cite{Pr}. 

Soit $\Lambda$ l'anneau des fonctions sym\'etriques, dont les sommes de puissances 
${\psi^{i}, i \ge 1}$ forment un syst\`eme de g\'en\'erateurs alg\'ebriques. 
On d\'efinit une action de $\Lambda$ sur tout polyn\^ome en 
posant
\[ \psi^{i} [\sum_{\alpha} c_{\alpha}u_{\alpha}]=\sum_{\alpha} 
c_{\alpha}u_{\alpha}^{i},\]
avec $c_{\alpha}$ constante r\'eelle et $u_{\alpha}$ un mon\^ome en 
$(v_{1},v_{2},\ldots,v_{n},y,t,\ldots)$.

Cette action s'\'etend naturellement \`a tout \'el\'ement de $\Lambda$. 
Pour toute partition $\mu$ et tout polyn\^ome $P$ on pose 
\[\psi^{\mu}[P]= \prod_{i=1}^{l(\mu)}\psi^{i}[P]=\prod_{i \ge 
1}{\psi^{i}[P]}^{m_{i}(\mu)}.\]
On a imm\'ediatement $\psi^{\mu}[PQ]=\psi^{\mu}[P]\psi^{\mu}[Q]$.

Si on note $\sigma^{n}$ la fonction sym\'etrique compl\`ete d'ordre $n$, 
on a pour tout  
polyn\^ome $P$,  la ``formule de Cauchy''
\[\sum_{n \ge 0}\sigma^{n}[P] = \sum_{\mu }\frac{\psi^{\mu}[P]}{z_{\mu}},\]
ou encore
\[\sigma^{n}[P] = \sum_{\left|{\mu }\right| = n}\frac{\psi^{\mu}[P]}{z_{\mu}}.\]

Nous avons pr\'ec\'edemment rencontr\'e cette formule g\'en\'erale 
sous la forme des deux ``sp\'ecialisations'' suivantes.

a) Soient une \textit{constante} $x$ et une \textit{ind\'etermin\'ee} $t$. La formule 
de Cauchy se formule comme suit
\[
\sum_{n \ge 0}\sigma^{n}[xt] =\sum_{m_i} \prod_{i} \frac{\displaystyle (xt^i)^{m_i}}
{\displaystyle {i}^{{m}_{i}}{m}_{i} !}=
\sum_{m_i} \prod_{i}  \frac {1}{m_i !} {\left( 
\frac{x t^i}{i}\right)}^{m_i}= (1-t)^{-x}.
\]
En d'autres termes, pour toute constante $x$ on a 
\[\sigma^{n}[x]=\binom{x+n-1}{n}.\]
 
b) Soient une \textit{constante} $x$ et deux \textit{ind\'etermin\'ees} $y$ et 
$Q$. On pose $q=Q-1$ et $P=xyq$. On a $\psi^{i}[xyq]=x\,y^i\psi^{i}[q]$ 
et $\psi^{i}[q]=\psi^{i}[Q-1]=Q^{i}-1.$ La formule de Cauchy se formule 
comme suit
\begin{multline*}
\sum_{n \ge 0}\sigma^{n}[xyq] =\sum_{m_i} {x}^{\sum_{i} {m}_{i}}\,{y}^{\sum_{i} 
i{m}_{i}}\prod_{i} \frac{\displaystyle {((1 + q)}^{i} - 1)^{m_i}}
{\displaystyle {i}^{{m}_{i}}{m}_{i} !}=\\
\sum_{m_i} \prod_{i}  \frac {1}{m_i !} {\bigg( x\,y^i\,
\frac{{(1 + q)}^{i} - 1}{i}\bigg)}^{m_i}.
\end{multline*}
Ce qui revient \`a dire que $\sigma^{n}[xq]$ est le coefficient de $y^n$ dans
\begin{equation*}
exp\left(x\,\sum_{i\ge 1} y^i\,\frac{{(1 + q)}^{i} - 1}{i}\right)
={\left(\frac{1-y}{1-y(1+q)}\right)}^{x}.
\end{equation*}
On retrouve ainsi la fonction apparue au Th\'eor\`eme 9.

De mani\`ere analogue, soient $z$ et $u$ deux ind\'etermin\'ees. On consid\`ere 
la formule de Cauchy
\begin{multline*}
\sum_{n \ge 0}\sigma^{n}[xyuq] =\sum_{m_i} {x}^{\sum_{i} {m}_{i}}\,{y}^{\sum_{i} 
i{m}_{i}} \prod_{i} {\psi^{i}[u]}^{m_{i}}\frac{\displaystyle {((1 + q)}^{i} - 1)^{m_i}}
{\displaystyle {i}^{{m}_{i}}{m}_{i} !}=\\
\sum_{m_i} \prod_{i}  \frac {1}{m_i !} {\bigg( x\,y^i\,
\frac{{(1 + q)}^{i} - 1}{i}\psi^{i}[u]\bigg)}^{m_i}.
\end{multline*}  
La fonction g\'en\'eratrice utilis\'ee \`a la 
Section 6 n'est autre que l'image du second membre par l'op\'erateur 
$\sum_{i} z^{i-1}
\psi^{i-1}[u]{d \over d\psi^{i}[u]}$, sp\'ecialis\'ee en $u=1$. 
Elle s'\'ecrit
\[x\,\left( \sum_{i\ge 1} z^{i-1} \frac{\psi^{i}[yq]}{i} \right)
\sum_{n \ge 0}\sigma^{n}[xyq].\]

\section{Motivations et perspectives}

Soient $\alpha$ un nombre r\'eel positif  et $(x_1,x_2,\ldots,x_N)$
une famille d'ind\'etermin\'ees ind\'ependantes. On dit qu'un 
polyn\^ome est ``sym\'etrique 
d\'ecal\'e'' s'il est sym\'etrique en les variables ``d\'ecal\'ees'' $x_{i}-i/\alpha$.
 Les ``polyn\^omes de Jack d\'ecal\'es''  \cite{La4} forment une base de l'alg\`ebre 
ainsi d\'efinie.
 
Chaque polyn\^ome  sym\'etrique d\'ecal\'e $f$ est 
enti\`erement d\'etermin\'e par ses valeurs $f(\lambda)$ pour les 
partitions $\lambda$.
Dans  \cite{La1,La3,La4} nous avons introduit la famille suivante de 
polyn\^omes sym\'etriques d\'ecal\'es, d\'efinie en la 
sp\'ecifiant sur les partitions. Pour 
toute partition $\lambda$ et tout entier $k\ge 0$, on note
\[d_{k}(\lambda) = \sum_{(i,j) \in \lambda} 
{\left(j-1-\frac{i-1}{\alpha}\right)}^k.\]
Pour tous entiers $j,k$ on pose
\[F_{jk}(\lambda)=\sum_{\left|{\mu }\right| = j} 
\frac{\displaystyle{\genfrac{\langle}{\rangle}{0pt}{}{\mu}{k}}}{z_{\mu}}\prod_{i \ge 
1}{d_{i}(\lambda)}^{m_{i}(\mu)}.\]

On d\'efinit ainsi des polyn\^omes sym\'etriques d\'ecal\'es  qui interviennent 
pour expliciter les ``polyn\^omes de Jack d\'ecal\'es'' 
(voir \cite{La4}, p.152--158). Dans 
\cite{La1} nous avons formul\'e \`a leur sujet plusieurs conjectures combinatoires. Il 
apparait que ces conjectures sont en fait \textit{ind\'ependantes} des 
quantit\'es $d_{k}(\lambda)$ et demeurent vraies si $d_{k}(\lambda)$ est 
remplac\'e par une ind\'etermin\'ee quelconque $X_{k}$.

Nous consid\`erons donc d\'esormais une famille (infinie) d'ind\'etermin\'ees ind\'ependantes
$X=(X_1,X_2,X_3,\ldots)$. Pour tous entiers $j,k \ge 0$ nous posons
\[P_{jk}(X)=\sum_{\left|{\mu }\right| = j} 
\frac{\displaystyle{\genfrac{\langle}{\rangle}{0pt}{}{\mu}{k}}}{z_{\mu}}\prod_{i \ge 
1}{X_{i}}^{m_{i}(\mu)}.\]

Comme on a $\genfrac{\langle}{\rangle}{0pt}{}{\mu}{k}  = 0$ si  
$ k < l(\mu )$, la sommation est limit\'ee aux partitions $\mu$ telles 
que $l(\mu) \le k$. Il en r\'esulte que $P_{jk}(X)$ est un polyn\^ome de degr\'e $k$. 
Comme on a $\genfrac{\langle}{\rangle}{0pt}{}{\mu}{k}  = 0$ si  
$k > \left|{\mu }\right|$, on a $P_{jk}(X)=0$ pour tout $k > j$. On pose par 
convention $P_{00}(X)=1$.

On a par exemple facilement 
\[P_{j1}(X)= X_{j},\]
\[P_{j2}(X)=\frac{1}{2}(j-1)X_{j}+\frac{1}{2}
\sum_{\begin{subarray}{1}j_{1}+j_{2}=j\\j_{1},j_{2} \ge 1\end{subarray}}X_{j_{1}}X_{j_{2}}.\]

Parmi les conjectures que nous avons formul\'ees  dans \cite{La1}, la conjecture 
centrale est la suivante (Conjecture 4). Elle explicite un 
d\'eveloppement en s\'erie formelle.
\begin{conj}
Soient $X_{0}$, $u$ et $X=(X_1,X_2,X_3,\ldots)$
des ind\'etermin\'ees ind\'ependantes. Pour tout entier $n \ge 1$ on a
\begin{multline*}
\sum_{\left|{\mu }\right| = n} \frac{(-1)^{n-l(\mu)}}{z_{\mu}}\prod_{i \ge 1}
 {\left(X_{0}+\sum_{p \ge 1}u^p\frac{{(i)}_{p}}{p!}\,X_{p} 
 \right)}^{m_{i}(\mu)} =\\
\sum_{j \ge 0} u^j\left( \sum_{k=0}^{min(n,j)} \binom{X_{0}-j}{n-k} P_{jk}(X) 
\right).
\end{multline*}
\end{conj}

Au membre de gauche il est \'evident que chaque terme correspondant 
\`a la partition $\mu$ est un polyn\^ome homog\`ene de degr\'e 
$l(\mu)$. On a en particulier
\begin{multline*}
\prod_{i \ge 1}{\left(X_{0}+\sum_{p \ge 1}u^p\frac{{(i)}_{p}}{p!}X_{p} 
 \right)}^{m_{i}(\mu)}=\\ 
 {X_{0}}^{l(\mu)}+X_{0}^{l(\mu)-1}\,\sum_{s \ge 1} u^s\,X_{s}\,
 \left( \sum_{i \ge 1} 
 m_{i}(\mu)\,\frac{{(i)}_{s}}{s!}\, \right)+\, 
 \textrm{autres termes}\ldots
\end{multline*}
 
Au membre de droite de la Conjecture 1, le terme en $X_{0}$  correspond au choix $j=k=0$. 
 On en d\'eduit l'\'egalit\'e
\[\sum_{\left|{\mu }\right| = n} \frac{{(-1)}^{n-l(\mu)}}{z_{\mu }} X_{0}^{l(\mu)}
= \binom{X_{0}}{n}.\]
On retrouve ainsi la relation classique cit\'ee au d\'ebut de la section 2. 
 
Au membre de droite de la Conjecture 1, le terme en  $X_{s}$ 
ne peut \^etre obtenu que pour $j=s$, et en restreignant la sommation dans 
$P_{sk}$ \`a la partition-ligne $(s)$. On a alors
\[P_{sk}(X)=\frac {\displaystyle{\binom{s}{k}}}{s} X_{s}+\, 
 \textrm{autres termes}\ldots\]
On en d\'eduit l'\'egalit\'e
\[
\sum_{\left|{\mu }\right| = n} {\frac{{(-1)}^{n-l(\mu)}}{{z}_{\mu }}}
\,X_{0}^{l(\mu)-1} \left( \sum_{i \ge 1} 
 m_{i}(\mu)\,\frac{{(i)}_{s}}{s!} \right) =
 \sum_{k = 1}^{min(n,s)} \binom{X_{0}-s}{n-k} \frac 
 {\displaystyle{\binom{s}{k}}}{s}.\]
On retrouve ainsi l'identit\'e d\'emontr\'ee au Th\'eor\`eme 4. 

Ces remarques nous conduisent \`a g\'en\'eraliser la Conjecture 1 sous 
la forme suivante. Cette nouvelle conjecture 
est triviale pour $r > n$ et on 
retrouve la Conjecture 1 pour $r=n$.
\begin{conj}
Soient $X_{0}$, $u$ et $X=(X_1,X_2,X_3,\ldots)$
des ind\'etermin\'ees ind\'ependantes. Pour tous entiers $n,r \ge 1$ on a
\begin{multline*}
\sum_{\left|{\mu }\right| = n} (-1)^{r-l(\mu)}
\frac{\displaystyle{\genfrac{\langle}{\rangle}{0pt}{}{\mu}{r}}}{z_{\mu}}
\prod_{i \ge 1}
 {\left(X_{0}+\sum_{p \ge 1}u^p\frac{{(i)}_{p}}{p!}\,X_{p} 
 \right)}^{m_{i}(\mu)} =\\
\sum_{j \ge 0} u^j \binom{n+j-1}{n-r}\left( \sum_{k=0}^{min(r,j)} 
\binom{X_{0}-j}{r-k} P_{jk}(X) 
\right).
\end{multline*}
\end{conj}

Comme pr\'ec\'edemment, l'\'egalit\'e des 
termes en $X_{0}$ s'\'ecrit
\[\sum_{\left|{\mu }\right| = n} {(-1)}^{r-l(\mu)} 
\frac{\displaystyle{\genfrac{\langle}{\rangle}{0pt}{}{\mu}{r}}}{z_{\mu }} X_{0}^{l(\mu)}
= \binom{n-1}{n-r} \binom{X_{0}}{r}.\]
On retrouve ainsi la relation (2). 
L'\'egalit\'e des termes en $X_{s}$ s'\'ecrit
\begin{multline*}
\sum_{\left|{\mu }\right| = n} {(-1)}^{r-l(\mu)} 
\frac{\displaystyle{\genfrac{\langle}{\rangle}{0pt}{}{\mu}{r}}}{z_{\mu }}
\,X_{0}^{l(\mu)-1} \left( \sum_{i \ge 1} 
 m_{i}(\mu)\,\frac{{(i)}_{s}}{s!} \right) = \\ \binom{n+s-1}{n-r}
 \sum_{k = 1}^{min(r,s)} \binom{X_{0}-s}{r-k} \frac 
 {\displaystyle{\binom{s}{k}}}{s}.
\end{multline*}
On retrouve ainsi l'identit\'e du Th\'eor\`eme 1.

La Conjecture 2 est triviale pour $r=1$. La sommation au membre de 
gauche est alors restreinte \`a la partition $\mu=(n)$, et la conjecture s'\'ecrit
\begin{multline*}
X_{0}+\sum_{p \ge 1}u^p\frac{{(n)}_{p}}{p!}\,X_{p} =
\sum_{j \ge 0} u^j \binom{n+j-1}{n-1} \sum_{k=0,1} 
\binom{X_{0}-j}{1-k} P_{jk}(X) \\=
X_{0}+ \sum_{j \ge 1} u^j \binom{n+j-1}{n-1} X_{j}.
\end{multline*}

Le lecteur pourra \'egalement v\'erifier la Conjecture 2 pour $r=2$. 
La sommation au membre de 
gauche est alors restreinte aux partitions de longueur $\le 2$, et la 
conjecture s'\'ecrit
\begin{multline*}
-\frac{n-1}{2}\left(X_{0}+\sum_{p \ge 1}u^p\frac{{(n)}_{p}}{p!}\,X_{p} 
 \right)+ \\
 \frac{1}{2}\sum_{i=1}^{n-1}\left(X_{0}+\sum_{p \ge 1}u^p\frac{{(i)}_{p}}{p!}\,X_{p} 
 \right)\left(X_{0}+\sum_{p \ge 1}u^p\frac{{(n-i)}_{p}}{p!}\,X_{p} 
 \right) = \\
 (n-1)\binom{X_{0}}{2}+\sum_{j \ge 1} u^j \binom{n+j-1}{n-2} \big( 
(X_{0}-j) X_{j}+ P_{j2}(X)\big).
\end{multline*}

Alain Lascoux a tout r\'ecemment obtenu une \'el\'egante 
d\'emonstration de la Conjecture 1 en utilisant les 
m\'ethodes de la th\'eorie des $\lambda$ -anneaux 
esquiss\'ees \`a la section 8. Il est vraisemblable que la 
Conjecture 2 pourra \^etre d\'emontr\'ee de la m\^eme mani\`ere \cite{La5}.

\end{document}